\documentclass[12pt]{amsart}
\usepackage{amsmath,amssymb,latexsym,graphicx}

\numberwithin{equation}{section}

\newcommand{\bd}{\mathbb{D}}

\newcommand{\br}{\mathbb{R}}

\newcommand{\bt}{\mathbb{T}}

\newcommand{\bz}{\mathbb{Z}}

\renewcommand{\a}{\alpha}

\renewcommand{\b}{\beta}

\renewcommand{\l}{\lambda}

\newcommand{\s}{\sigma}

\renewcommand{\r}{\rho}

\newcommand{\p}{\varphi}

\newcommand{\pp}{\Phi}

\newcommand{\pps}{\Psi}
\renewcommand{\t}{\tau}

\renewcommand{\th}{\theta}

\renewcommand{\d}{\delta}

\newcommand{\oo}{\Omega}

\newcommand{\g}{\gamma}

\newcommand{\ep}{\varepsilon}

\newcommand{\z}{\zeta}

\renewcommand{\k}{\kappa}

\newcommand{\nt}{\noindent}

\newcommand{\bsl}{\backslash}

\newcommand{\ovl}{\overline}

\newcommand{\pt}{\partial}

\newcommand{\ti}{\tilde}

\newcommand{\lt}{\left}

\newcommand{\rt}{\right}

\newcommand{\dsp}{\displaystyle}

\renewcommand{\ss}{\mathrm{(S)}}

\newcommand{\ps}{\mathrm{(pS)}}
\newcommand{\cf}{\mathcal{F}}

\newcommand{\cc}{\mathcal{C}}

\newcommand{\cp}{\mathcal{P}}

\newcommand{\tr}{\mathrm{tr}\,}

\newcommand{\re}{\mathrm{Re}\,}

\newcommand{\im}{\mathrm{Im}\,}

\newcommand{\rk}{\mathrm{rank}\,}

\renewcommand{\sl}{\mathrm{s-lim}\,}

\newcommand{\epr}{{\hfill $\Box$}}

\newtheorem{theorem}{Theorem}[section]
\newtheorem{lemma}[theorem]{Lemma}

\newtheorem{remark}[theorem]{Remark}

\begin{document}
\title[Asymptotics of orthogonal polynomials]
{Asymptotics of the orthogonal polynomials\\
for the Szeg\H{o} class with a polynomial weight}

\author{S. Denisov}
\author{S. Kupin}

\email{denissov@its.caltech.edu}
\email{kupin@cmi.univ-mrs.fr}


\date{April 22, 2004}

\address{Department of Mathematics 253-37, Caltech, Pasadena CA 91125, USA.}
\address{CMI, Universit\'e de Provence, 39 rue Joliot  Curie, 13453
  Marseille Cedex 13, France.}

\begin{abstract}
Let $p$ be a trigonometric polynomial, non-negative on the unit
circle $\mathbb{T}$. We say that a measure $\sigma$ on
$\mathbb{T}$ belongs to the polynomial Szeg\H{o} class, if
$d\sigma(e^{i\theta})=\sigma_{ac}'(e^{i\theta})d\theta+d\sigma_s(e^{i\theta})$,
$\sigma_s$ is singular, and

$$
\int^{2\pi}_0 p(e^{i\th})\log \s'_{ac}(e^{i\th})\, d\th>-\infty
$$
For the associated orthogonal polynomials $\{\varphi_n\}$, we
obtain pointwise asymptotics inside the unit disc $\mathbb{D}$.
Then we show that these asymptotics hold in $L^2$-sense on the
unit circle. As a corollary, we get an existence of certain
modified wave operators.
\end{abstract}

\maketitle

\vspace{-0.7cm}
\section*{Introduction}
\label{s0}


Let $\s$ be a non-trivial Borel probability measure on  the unit
circle $\bt=\{z:|z|=1\}$. Consider polynomials $\{\p_n\}$
orthonormal with respect to $\sigma$,
$$
\int_\bt\p_n\ovl{\p_m}\, d\s=\d_{nm}
$$
where $\d_{nm}$ is the Kronecker's symbol. Sometimes, it is more convenient to
work with monic orthogonal polynomials $\{\pp_n\},\
\pp_n(z)=z^n+a_{n,n-1}z^{n-1}+\ldots+a_{n,0}$. These polynomials satisfy
$$
\int_\bt\pp_n\ovl{\pp_m}\, d\s=c_n\d_{nm}
$$
with $c_n=||\pp_n||^2_\s=\int_\bt |\pp_n|^2\, d\s$.

It is well-known \cite{ge1,si1} that polynomials $\{\pp_n\}$
generate a sequence \mbox{$\{\a_n\}, |\a_n|<1$},  of the so-called
Verblunsky coefficients through the recurrence relations
$$
\left\{
\begin{array}{lcl}
\pp_{n+1}(z)&=&z\pp_n(z)-\bar\a_n\pp^*_n(z)\\
\pp^*_{n+1}(z)&=&\pp^*_n(z)-\a_nz\pp_n(z)
\end{array}
\right.
$$
where
$\pp_0(z)=1, \pp^*_0(z)=1$, and $\pp^*_n(z)=z^n\ovl{\pp_n(1/\bar z)}$.
Conversely, the measure $\s$ (and polynomials $\{\p_n\}$) are completely
determined by the sequence $\{\a_k\}$ of its Verblunsky
parameters. Hence, it is natural to study properties of the
sequence $\{\a_k\}$ and polynomials $\{\p_n\}$ in terms of $\s$
and vice versa.

We say that $\s$ is a Szeg\H o measure ($\s\in\ss$, for brevity),
if $d\s=d\s_{ac}+d\s_s=\s'_{ac}dm+d\s_s$ and the density
$\s'_{ac}$ of the absolutely continuous part of $\s$ is such  that
$$
\int_\bt \log \s'_{ac}\, dm>-\infty
$$
Here, the singular part of $\s$ is denoted by $\s_s$, and $m$ is the
probability Lebesgue  measure on $\bt$, $dm(t)=dt/(2\pi it)=1/(2\pi)\,
d\th, t=e^{i\th} \in\bt$.

The following theorem is classical.

\begin{theorem}[{\cite{ge1,sz}}]\label{t01}
The following assertions are equivalent
\begin{itemize}
\item[\it i)\ \ ] the sequence $\{\a_k\}$ is in
$\ell^2(\mathbb{Z}_+)$, \item[\it ii)\ ] the measure $\s$ belongs
to the Szeg\H o class, \item[\it iii)] analytic polynomials are
not dense in $L^2(\s)$.
\end{itemize}
\end{theorem}

We denote by $\cp_0$ the set of analytic polynomials $f$ such that
$f\not=0$ on $\bd$ and $f(0)>0$. Let also $\cp_1=\{f: f\in\cp_0,
f(0)=1\}$. Then, the last statement of the theorem can be made
more precise. Namely, we have \cite{ge1,sz} that
\begin{eqnarray}
d(\cp_1,0)^2_{L^2(\s)}&=&\inf_{f\in \cp_1} ||f||^2_\s=
\inf_{\dsp f\in\cp_0, ||f||_\s\le1} |f(0)|^{-2} \label{e02}\\
&=&\exp \int_\bt\log\s'_{ac}\, dm \nonumber
\end{eqnarray}

If $\s\in\ss$, we  define a function $D$, lying in the Hardy space $H^2(\bd)$ on the
unit disk $\bd=\{z:|z|<1\}$, as
\begin{equation}\label{e01}
D(z)=\exp\lt(\frac12 \int_\bt \frac{t+z}{t-z}\, \log\s'_{ac}(t)\,dm(t)\rt)
\end{equation}

\begin{theorem}[{\cite{ge1,sz}}]\label{t02}
Let $\s\in\ss$. Then
$$
\lim_{n\to\infty} D(z)\p^*_n(z)=1
$$
for every $z\in\bd$. Moreover,
$$
\lim_{n\to\infty} \int_\bt |D\p^*_n-1|^2\, dm=0
$$
\end{theorem}

A modern presentation and recent advances in this direction can be found in
\cite{khr1,si1}.

It seems interesting to obtain similar results for
different classes of measures. Consider
a trigonometric polynomial $p$ with the property $p(t)\ge 0,
t\in\bt$. Without loss of generality we can assume it is in the form
\begin{equation}\label{e03}
p(t)=\prod_{k=1}^N |t-\z_k|^{2\kappa_k}
\end{equation}
where $\{\z_k\}$ are points on $\bt$ and $\kappa_k$ are their
``multiplicities''. We say that $\s$ is in the polynomial Szeg\H o
class (i.e., $\s$ is a (pS)-measure or $\s\in\ps$), if
$d\s=\s'_{ac}dm+d\s_s$, $\s_s$ being the singular part of the
measure, and
\begin{equation}\label{e04}
\int_\bt p(t)\log\s'_{ac}(t)\, dm(t)>-\infty
\end{equation}

The main result of the paper is a counterpart of Theorem \ref{t02}
for orthogonal polynomials with respect to polynomial Szeg\H o
measures. We want to mention here that similar results for Jacobi
matrices, including $L^2$-asymptotics of orthogonal polynomials on
a segment, were obtained recently by Damanik and Simon
\cite{ds}. Authors were able to deal with the case considered in
\cite{dks}, Theorem 3.1.

We introduce some notations. Actually, all objects appearing below
should be indexed by the polynomial $p$ from \eqref{e03}. We omit
this dependence.

Let $\s\in\ps$. Consider a modified Schwarz kernel
\begin{equation}\label{e041}
K(t,z)=\frac{t+z}{t-z}\,\frac{q(t)}{q(z)}=\frac{t+z}{t-z}\,\frac{q_0(t)}{q_0(z)}
\end{equation}
where $q_0(t)=\prod^N_{k=1}(t-\z_k)^{2\k_k}/t^{N'},
N'=\sum_k\k_k$, and $q(t)=Cq_0(t)$. The constant
$C$ equals $(\prod_k(-\z_k)^{\k_k})^{-1}$, so that $|C|=1$ and
$q(t)=\prod_k|t-\z_k|^{2\k_k}=p(t)$ for $t\in\bt$. Let
us introduce
\begin{eqnarray}
\ti D(z)&=&\exp\lt(\frac12\int_\bt K(t,z)\log\s'_{ac}(t)\, dm(t)\rt)
\label{e042}\\
\ti\p^*_n(z)&=&\exp\lt(\int_\bt K(t,z)\log|\p^*_n(t)|\, dm(t)\rt)
\label{e043}
\end{eqnarray}
with $z\in\bd$. We call functions $\{\ti\p^*_n\}$ the modified
reversed polynomials. The properties of the kernel $K$ easily
imply that $|\ti D|^2=\s'_{ac}$ and $|\ti\p^*_n|=|\p^*_n|$ a.e.~on
$\bt$, see Lemma \ref{l3}. It is also useful to consider the
functions
$$
\psi_n(z)=\frac{\ti\p^*_n(z)}{\p^*_n(z)}=
\exp\lt(\int_\bt\frac{t+z}{t-z}\lt(\frac{q(t)}{q(z)}-1\rt)\log |\p^*_n(t)|\,
dm(t)\rt)
$$
Clearly $|\psi_n|=1$ a.e.~on $\bt$ and, similarly to $ii)$, Lemma
\ref{l3}
\begin{equation}\label{e05}
\psi_n(z)=\exp\lt\{A_{0n}+\sum^N_{k=1}\sum^{2\k_n}_{j=1}
A_{j,kn}\lt(\frac{z+\z_k}{z-\z_k}\rt)^j\rt\}
\end{equation}
where $A_{0n},A_{2j, kn}\in i\br$ and $A_{2j+1,kn}\in\br$. The coefficients
$\{A_{0n}, A_{j,kn}\}$ can be expressed in a closed form through the Verblunsky
coefficients $\{\a_k\}$.

The following theorem holds.

\begin{theorem}\label{t2} Let $\s\in\ps$. Then
$$
\lim_{n\to\infty} \ti D(z)\ti\p^*_n(z)=1
$$
for every $z\in\bd$.
\end{theorem}

The proof of the theorem is largely inspired by the classical
proof of Theorem \ref{t02} \cite{ge1,sz}, and it is based on
appropriate sum rules. These sum rules  are obtained in Theorem
\ref{t1}. Their proof is a translation of \cite{nvyu}, Theorem 1.4
to the case of orthogonal polynomials on the unit circle. We also
mention that the relations we prove in the theorem are closely
related to sum rules obtained in \cite{ks1,ku1,ku2}. A counterpart
of Theorem \ref{t2} for Jacobi matrices is \cite{nvyu}, Theorem
1.5.

A subsequent analysis shows that Theorem \ref{t2} can be considerably
strengthened.

\begin{theorem}\label{t3} Let $\s\in\ps$. Then
$$
\lim_{n\to\infty} \int_\bt |\ti D\ti\p^*_n-1|^2\, dm=0
$$

\end{theorem}

The proof of the theorem is rather technical. One of the main
observations leading to the statement is that
$$
\lim_{n\to\infty}\int_I|\ti D\ti\p^*_n-1|^2\, dm=0
$$
for any closed arc $I\subset \bt$ that does not contain points $\{\z_k\}$.
We prove the latter relation showing that
$$
|\ti D\ti\p^*_n(z)|\le\frac{C_\ep}{\sqrt{1-|z|}}
$$
for $z\in \bd\bsl(\cup_k B_\ep(\z_k))$, $B_\ep(\z)=\{z:
|z-\z|<\ep\}$, whenever $\ep>0$ is small enough. It is crucial that the above
constant $C_\ep$ does not depend on $n$.

We apply Theorem \ref{t3} to construct modified wave operators for
the CMV-represen\-tations $\cc,\cc_0$ associated to measures
$\s\in\ps$ and $m$, see Section \ref{s11} for the definitions and
notation. For the Szeg\H{o} case, 
the classical wave operators were described recently by
Simon \cite[Ch.~10]{si1}.  Let $\cf_0:L^2(m)\to
\ell^2(\bz_+), \cf:L^2(\s)\to \ell^2(\bz_+)$ be the Fourier
transforms related to $\cc$ and $\cc_0$.  Recall that
$$
\cc=\cf z\cf^{-1},\quad \cc_0=\cf_0 z\cf^{-1}_0
$$

\begin{theorem}\label{t4}
Let $\s\in\ps$. The limits
\begin{equation}\label{e06}
\ti\oo_\pm=\sl_{n\to\pm\infty}\, e^{W(\cc,2n)}\cc^n\cc_0^{-n}
\end{equation}
exist. Here
$$
W(\cc,n)=A_{0n}+\sum^N_{k=1}\sum^{2\k_k}_{j=1}
A_{j,kn}\lt(\frac{\cc+\z_k}{\cc-\z_k}\rt)^j
$$
and coefficients $\{A_{0n},A_{j,kn}\}$ are defined in
\eqref{e05}. Furthermore,
\begin{equation}\label{e141}
\cf^{-1}\ti\oo_+\cf_0=\chi_{E_{ac}}\,\frac1{\ti D},\quad
\cf^{-1}\ti\oo_-\cf_0=\chi_{E_{ac}}\,\frac1{\bar{\ti D}}
\end{equation}
where $E_{ac}=\bt\bsl\mathrm{supp}\, \s_s$.
\end{theorem}

In the formulation above, $\sl$ refers to the limit in
$\ell^2(\bz_+)$ in the strong sense. A natural problem is to pass
from wave operators \eqref{e06} to operators of the form
$$
\sl_{n\to\pm\infty}\, \cc^n\cc_0^{-n}e^{\ti W(\cc_0,n)}
$$
This question is still open, see \cite{chki} in this connection.

Finally, we address a variational principle that is naturally connected to
measures from a (pS)-class. Let $p$ be the trigonometric polynomial from
\eqref{e03}. We pick a constant $C_0$ in a way that
$C_0\int_\bt p\, dm=1$, and let $p_0=C_0p$.

For a $g\in \cp_0$, we define
$$
\l(g)=\exp\lt(\int_\bt p_0\log|g|\, dm\rt)
$$
and $\cp'_1=\{g: g\in\cp_0,\ \l(g)=1\}$.

\begin{theorem}\label{t5}
Let $d\s=\sigma_{ac}'\, dm+d\s_s$. Then
\begin{eqnarray}
\exp\lt(\int_\bt p_0\log\frac
{\sigma_{ac}'}{p_0}\,dm\rt)&\le&\inf_{g\in\cp'_1}||g||^2_\s=
\inf_{ \dsp
\begin{array}{c}
g\in\cp_0,\\
||g||_\s\le 1
\end{array}}
\frac1{|\l(g)|^2} \label{e15}\\
&\le&\exp\lt(\int_\bt p_0\log \sigma_{ac}'\, dm\rt) \nonumber
\end{eqnarray}
\end{theorem}

Remind that $\s$ is a Szeg\H o measure if and only if the
system $\{e^{iks}\}_{k\in\bz}$ is uniformly minimal in
$L^2(\s)$ \cite[Ch.~3]{ga}, \cite[Ch.~6]{nk}. Saying that $\s$ is a (pS)-measure
translates into the
uniform minimality of another system, $\{e^{ik\nu(s)}\}_{k\in\bz}$,
in the same
space $L^2(\s)$. Above,
$$
\nu(s)=\int^s_0 p_0(e^{is'})\, ds'
$$
where $s,s'\in [0,2\pi]$; see \cite{nvyu}, Lemma 2.2.

We now turn to the concrete example to illustrate our results. It
was proved recently in \cite[Ch. 2]{si1} that
$\s\in\mathrm{(p}_1\mathrm{S)}$ with
$$
p_1(t)=\frac12\, |1-t|^2=1-\cos\th
$$
if and only if $\{\a_k\}\in \ell^4(\mathbb{Z}_+)$ and
$\{\a_{k+1}-\a_k\}\in \ell^2(\mathbb{Z}_+)$ (above, $t=e^{i\th}$).
This class of parameters was studied earlier in \cite{de1}.
Theorems \ref{t2}--\ref{t5} readily apply to this special case. In
particular, we have
\begin{eqnarray*}
K_1(t,z)&=&\frac{t+z}{t-z}\, \frac{(t-1)^2}{t}\, \frac{z}{(1-z)^2}
=-\frac{z}{(1-z)^2}\frac{t+z}{t-z}\,|1-t|^2\\
\ti D_1(z)&=&\exp\lt(\frac12\int_\bt K_1(t,z)\log\s'_{ac}(t)\, dm(t)\rt)
\end{eqnarray*}
and
$$
\psi_n(z)=\exp\lt(A_n\frac{1+z}{1-z} +
B_n\lt\{\lt(\frac{1+z}{1-z}\rt)^2-1\rt\}\rt)
$$
where
$$
A_n=\sum^n_{k=0}\log (1-|\a_k|^2)^{1/2},\quad
B_n=\frac i4\, \im\Big(\a_0-\sum^n_{k=1}\bar\a_{k-1}\a_k\Big)
$$
Recently, for $t_1,t_2\in\bt$, the following class of polynomials $p$ was
considered \cite{sz1}
$$
p(t)=|(t-t_1)(t-t_2)|^2
$$
and the criteria for \eqref{e04} to be true were obtained in terms
of the Verblunsky coefficients. Methods of the current paper are
also applicable to this case.

We conjecture that counterparts of Theorems \ref{t3}, \ref{t4}
hold true for Jacobi matrices; see \cite{ks1,nvyu} in this
connection.

The paper is organized as follows. The preliminaries are in
Section \ref{s11}. The sum rules we use in the proof of Theorem
\ref{t2} are obtained in Section \ref{s1}. Theorem \ref{t2} itself
is proved in Section \ref{s2}, and it is ``upgraded'' to the
asymptotics in $L^2(\bt)$-sense in Section \ref{s3}. Section
\ref{s4} deals with the modified wave operators and the
variational principle from Theorem \ref{t5}.

As usual, $H^p(\bd)$ is the Hardy space of analytic functions on
the unit disk \cite{ga}. For an arc $I\subset \bt$, we write
$L^2(I)$ to refer to the standard $L^2$-space with the Lebesgue
measure on $I$. We set $\log^+ x=(|\log x|+\log x)/2$ and $\log^-
x=(|\log x|-\log x)/2$ for $x>0$. Also, $C$ is a constant changing
from one relation to another.

\medskip\nt{\it Acknowledgments.}\ \ The
authors are grateful to N.~Nikolskii, B.~Simon, and P.~Yuditskii
for helpful discussions.

\section{Preliminaries}
\label{s11} It is useful to keep in mind the simple general
properties of the measures from a (pS)-class. Following
\cite{khr1,si1}, we say that $\sigma$ belongs to the Erd\H{o}s
class ($\s\in\mathrm{(E)}$) if $\s'_{ac}>0$ a.e. on $\bt$. A measure is in the
Nevai class ($\s\in\mathrm{(N)}$) if $\lim_{n\to\infty}\alpha_n=0$. Lastly,
$\s$ is a Rakhmanov measure (i.e., $\s\in\mathrm{(R)}$) if
\begin{equation}\label{e104}
\mathrm{w-}\lim_{n\to\infty} |\p_n|^2\,d\s=dm
\end{equation}
The following relations are true \cite[Sect.~2, 6, 7]{khr1},
\cite[Ch.~7]{si1}
\begin{equation}
\mathrm{(S)}\subset\mathrm{(pS)}\subset\mathrm{(E)}\subset\mathrm{(N)}
\subset\mathrm{(R)}
\label{e103}
\end{equation}
Here, the first and second inclusions are obvious.

Let us recall a few facts on the so-called CMV-represen\-tations.
More information on the topic can be found in
\cite[Ch.~4]{cmv,si1}.

Let $\s$ be a measure on $\bt$. Consider the unitary operator
$U:L^2(\s)\to L^2(\s)$ given by the formula $Uf(t)=tf(t), f\in
L^2(\s)$. It turns out one can find an orthonormal basis
$\{\chi_n\}_{n\in\bz}$ in $L^2(\s)$ such that the matrix of $U$ in
this basis has a reasonably simple form. Namely, we set for
$n\in\bz_+=\{0,1,2,\ldots\}$
$$
\chi_n(z)=\lt\{
\begin{array}{lcl}
z^{-k}\p^*_{2k}(z),&& n=2k\\
z^{-(k-1)}\p_{2k-1}(z),&& n=2k-1
\end{array}
\right.
$$

\begin{theorem}[{\cite[Ch.~4]{cmv,si1}}]\label{t11}
The operator $U$, defined above, is unitarily equivalent to the
operator $\cc:\ell^2(\bz_+)\to \ell^2(\bz_+)$ of the form
$$
\cc=\cc(\s)=\lt[
\begin{array}{cccccc}
*&*&*&0&0&\ldots\\
*&*&*&0&0&\ldots\\
0&*&*&*&*&\ldots\\
0&*&*&*&*&\ldots\\
0&0&0&*&*&\ldots\\
\vdots&\vdots&\vdots&\vdots&\vdots&\ddots\\
\end{array}
\rt]
=
\lt[
\begin{array}{ccc}
A_0&0&\ldots\\
0&A_1&\ldots\\
\vdots&\vdots&\ddots
\end{array}
\rt]
$$
where $\a=\{\a_k\}$ is the sequence of Verblunsky coefficients of $\s$,
\begin{eqnarray*}
A_j&=&
\lt[
\begin{array}{cccc}
\bar\a_{k+1}\rho_k&-\bar\a_{k+1}\a_k&\bar\a_{k+2}\rho_{k+1}&\rho_{k+2}\rho_{k+1}
\\
\rho_{k+1}\rho_k&-\rho_{k+1}\a_k&-\bar\a_{k+2}\a_{k+1}&-\rho_{k+2}\a_{k+1}
\end{array}
\rt]\\
A_0&=&
\lt[
\begin{array}{ccc}
\bar\a_0&\bar\a_1\rho_0&\rho_1\rho_0\\
\rho_0&-\bar\a_1\a_0&-\rho_1\a_0
\end{array}
\rt]
\end{eqnarray*}
and $\rho_k=(1-|\a_k|^2)^{1/2}$.
\end{theorem}

The matrix $\cc$ is called a CMV-representation associated to the measure $\s$.

It is easy to see that the map $\cf: L^2(\s)\to \ell^2(\bz_+)$,
carrying out the unitary equivalence
$$
\cc=\cf z\cf^{-1}
$$
is determined by relations $(\cf f)_n=\int_\bt f\bar\chi_n\, d\s$,
where $f\in L^2(\s)$. Similar objects, related to the Lebesgue
measure $m$, are indexed by 0. That is, its CMV-matrix is denoted
by $\cc_0$, $\{\chi^{(0)}_n\}$ and $\cf_0$ are the standard basis
$\{t^k\}_{k\in\bz}$ and the Fourier transform, respectively.

Now, we denote by $\cc_n$ the $n\times n$ upper left block of
$\cc$. One can prove that \cite{si1}, Theorems 1.7.14, 4.2.14
$$
\p_n(z)=\frac1{A_n}\,\det(z-\cc_n)(z-\cc_{0,n})^{-1},\quad
\p^*_n(z)=\frac1{A_n}\,\det(1-z\bar\cc_n)(1-z\bar\cc_{0,n})^{-1}
$$
with $A_n=\prod^{n-1}_{k=0} \r_k$. Recalling definition \eqref{e01}, we get the
following theorem.

\begin{theorem}[{\cite{si1}, Theorem 4.2.14}]
\label{t12}
Let $\sum_k|\a_k|<\infty$. Then, for $z\in\bd$,
$$
D(z)=A_\infty\, \det(1-z\bar\cc_0)(1-z\bar\cc)^{-1}
$$
Moreover, we have $\log D(z)=t_0+\sum^\infty_{k=1}(t_k/k)\, z^k$
and
\begin{equation}\label{e07}
t_0=\sum_k\log \r_k=\sum_k\log(1-|\a_k|^2)^{1/2},\quad
t_k=\tr(\bar\cc^k-\bar\cc^k_0)
\end{equation}
with $k\ge1$.
\end{theorem}

\section{Polynomial Szeg\H o condition and corresponding sum rules}
\label{s1}

We fix the polynomial $p$ \eqref{e03} for the rest of this paper.
For the sake of transparency, we also assume $\kappa_k=1$; the
discussion of the general case follows the same lines.

The goal of this section is to obtain the sum rules similar to
\cite{si1}, Theorem 2.8.1 and \cite{nvyu}, Theorem 1.4. With the
exception of simple technical details, our argument follows
word-by-word a reasoning from \cite{nvyu}.

We start with a CMV-representation $\cc$ having the property
$\rk(\cc-\cc_0)<\infty$. Note that this is equivalent to saying that the
sequence of Verblunsky coefficients $\{\a_k\}$, corresponding to $\cc$, is
finite. Therefore, $\sum_k|\a_k|<\infty$ and, by Theorem \ref{t12},
$$
\log D(z)=t_0+\sum^\infty_{k=1} \frac{t_k}{k}\, z^k
$$
with coefficients $\{t_k\}$ given by \eqref{e07}. Since $\log D\in
H^\infty(\bd)\cap C(\bar\bd)$, this yields
\begin{equation}\label{e102}
\int_\bt \log|D|^2\, dm=2t_0,\quad \int_\bt t^k\log|D|^2\, dm=\frac{\bar t_k}k
\end{equation}
Taking polynomial $p$ from \eqref{e03}, we define an analytic
polynomial $P$ through the relations
\begin{equation}\label{e101}
p_1=2P_+(p),\quad P'(t)=\frac{p_1(t)-p_1(0)}{t},\quad P(0)=0
\end{equation}
here $P_+:L^2(\bt)\to H^2(\bd)$ is the Riesz projection \cite[Ch.~3]{ga}.
\begin{lemma}\label{l1}
Let $p$ be as above and $\rk(\cc-\cc_0)<\infty$. We have
\begin{equation}\label{e2}
\int_\bt p(t)\log|D(t)|^2\, dm(t)=A_0t_0+\re\tr(P(\cc)-P(\cc_0))
\end{equation}
where $A_0=p_1(0)=2\int_\bt p(t)\, dm$.
\end{lemma}

\begin{proof}
Write $p$ as $p(t)=a_0+2\re\sum_{j=1}^N a_jt^j$. Recalling
\eqref{e102}, we get
\begin{eqnarray*} \int_\bt p\log |D|^2\,
dm&=&2a_0t_0+2\re\sum_{j=1}^N a_j\int_\bt t^j\log|D|^2\,
dm\\
&=&2a_0t_0+2\re\sum_{j=1}^N\frac{a_j}{j}\bar t_j=2a_0t_0+2\re\sum_{j=1}^N
\frac{a_j}j\tr(\cc^j-\cc_0^j)
\end{eqnarray*}
It remains to notice that the polynomial 2$\sum_j (a_j/j) z^j$
above is indeed $P$ given by \eqref{e101} and $A_0=2a_0$. Hence,
the last expression in the displayed formula is exactly
$A_0t_0+\re\tr(P(\cc)-P(\cc_0))$, and the lemma is proved.
\end{proof}

Let us set
\begin{eqnarray*}
\pp(\cc)&=&\int_\bt p\log \s'_{ac}\, dm\\
\pps(\cc)&=&A_0t_0+\re\tr(P(\cc)-P(\cc_0))
\end{eqnarray*}
Notice that $\pp(\cc)$ is exactly the left-hand side of equality \eqref{e2}.

We now rewrite $\pps(\cc)$ in a different form. Since $t_0=\sum
\log\r_k$, we have
\begin{equation}\label{e3}
\pps(\cc)=\sum^\infty_{k=0}\lt\{A_0\log\r_k+\re ((P(\cc)-P(\cc_0))e_k,e_k)\rt\}
\end{equation}
here $\{e_k\}$ is the standard basis in $\ell^2(\bz_+)$. Consider
the shift $S:\ell^2(\bz_+)\to \ell^2(\bz_+)$ given by
$Se_k=e_{k+1}$. For a bounded operator $A$ on $\ell^2(\bz_+)$,
take $\t(A)=S^*AS$. It is obvious that the matrix of $\t^k(A)$ is
obtained from the matrix of $A$ by dropping its first $k$ rows and
columns.

Furthermore, the degree of the polynomial $P$  is $N$, the matrix $\cc$ is
five-diagonal, so $P(\cc)$ contains $4N+5$ non-zero diagonals. Consequently,
equality \eqref{e3} is exactly the same as
$$
\pps(\cc)=\sum_{k=0}^{2N+1}\{A_0\log\r_k+\re((P(\cc)-P(\cc_0)e_k,e_k)\}+
\sum^\infty_{k=0} \psi\circ\t^k(\cc)
$$
where
$$
\psi(\cc)=A_0\log\r_{2N+2}+((P(\cc)-P(\cc_0))e_{2N+2},e_{2N+2})
$$
is a function of a finite number of Verblunsky
coefficients.

The following lemma is similar to \cite{nvyu}, Lemma 3.1.
\begin{lemma}[{\cite{nvyu}}]\label{l2}
There exists a function $\g$ depending on $l=4N+4$ arguments such that
$$
\psi(x_1,\ldots,x_{l+1})=\eta(x_1,\ldots,x_{l+1})-\g(x_2,\ldots,x_{l+1})+\g(x_1,
\ldots,x_{l})
$$
and $\eta(x_1,\ldots,x_{l+1})\le 0$ for any collection $(x_1,\ldots,x_{l+1})$.
\end{lemma}
The proof of the lemma relies on the fact that $\pps(\cc)\le
C<\infty$ for all $\cc$ with the property $\rk(\cc-\cc_0)<\infty$.
This is obviously true because we have $\pps(\cc)=\pp(\cc)$ for
these $\cc$ and $\pp(\cc)$ is uniformly bounded away from $\infty$
by the Jensen inequality. Now, define
\begin{eqnarray}\label{e4}
\ti\pps(\cc)&=&\sum_{k=0}^{2N+1}\{A_0\log\r_k+\re((P(\cc)-P(\cc_0)e_k,e_k)\}\\
&+&\sum^\infty_{k=0} \eta\circ\t^k(\cc)+\g(\cc)\nonumber
\end{eqnarray}

\begin{theorem}[{\cite{nvyu}}]\label{t1}
A measure $\s$ lies in the (pS)-class (see \eqref{e04}) if and only if
$\ti\pps(\cc)>-\infty$, or, equivalently,
$\sum^\infty_{k=0}\eta\circ\t^k(\cc)>-\infty$. Moreover,
\begin{equation}\label{e41}
\pp(\cc)=\ti\pps(\cc)=\pps(\cc)
\end{equation}
\end{theorem}

The proof literally follows \cite{nvyu}, Theorem 1.4 and it is
close in spirit to arguments from \cite{ks1}, \cite{si1}, Theorem
2.8.1. Its main ingredients are the non-positivity of $\eta$ in
\eqref{e4} and  the fact that $\lim_{k\to\infty}\a_k=0$.

\section{Pointwise asymptotics for orthogonal polynomials\\ on the unit disk}
\label{s2}

We start with the following lemma.
\begin{lemma}\label{l3}
Let $\s\in\ps$, the polynomials $\ti\p^*_n$ and the function $\ti D$ be defined
in \eqref{e043}, \eqref{e042}. Then
\begin{itemize}
\item[\it i)\ ] $|\ti D(t)|^2=\s'_{ac}(t)$ a.e. on $\bt$,
\item[\it ii)]  $\ti\p^*_n=\psi_n\p^*_n$ and $|\psi_n(t)|=1$ a.e.
on $\bt$. Moreover,
\begin{equation}\label{e7}
\psi_n(z)=\exp\lt(A_{0n}+\sum_k\lt(A_{1,kn}\frac{z+\z_k}{z-\z_k}+
A_{2,kn}\lt\{\frac{z+\z_k}{z-\z_k}\rt\}^2\rt)\rt)
\end{equation}
where $A_{0n},A_{2,kn}\in i\br$ and $A_{1,kn}\in\br$.
\end{itemize}
\end{lemma}

\begin{proof}
To prove claim $i)$, observe that
$$
\log|\ti D(z)|^2=\re\int_\bt\frac{t+z}{t-z}\frac{q(t)}{q(z)}\log\s'_{ac}(t)\,
dm(t)
$$
Also, $\im\{q(t)/q(z)\}$ tends uniformly to $0$ as $z$ goes to
$\bt\bsl\{\z_k\}$, and $q=\re q=p$ on $\bt$.
Consequently, for a.e.~$t_0\in\bt\bsl\{\z_k\}$,
\begin{eqnarray*}
\lim_{z\to t_0}\log |\ti D(z)|^2&=&\lim_{z\to t_0} \frac1{\re q(z)}\int_\bt
\re\frac{t+z}{t-z}\, \re q(t)\log\s'_{ac}(t)\, dm(t)\\
&=&\frac1{p(t_0)}\, p(t_0)\log\s'_{ac}(t_0)
\end{eqnarray*}
where we used the standard properties of the Poisson kernel $\re
(t+z)/(t-z)$.

The computation also shows that $|\ti\p^*_n|=|\p^*_n|$ a.e. on
$\bt$, and, in particular, $|\psi_n|=1$ a.e. On the other hand,
$$
\psi_n(z)=\exp\lt(\int_\bt\lt\{\frac{t+z}{t-z}\,
\frac{q(t)-q(z)}{q(z)}\rt\}\log|\p^*_n(t)|\, dm(t)\rt)
$$
The function in the curled brackets is rational with respect to
$z$ and its degree is $2N$. Its poles have multiplicities two and
they are located at $\{\z_k\}$. So, we get
$$
\frac{t+z}{t-z}\, \frac{q(t)-q(z)}{q(z)}=a_0(t)+\sum_k
\lt(a_{1k}(t)\frac{z+\z_k}{z-\z_k}+a_{2k}(t)\lt\{\frac{z+\z_k}{z-\z_k}\rt\}^2\rt
)$$ where $a_0,a_{1k},a_{2k}$ are some trigonometric polynomials
(i.e., polynomials with respect to $t,\bar t$). We now put
$A_{0n}=\int_\bt a_0\log|\p^*_n|\, dm, A_{1,kn}=\int_\bt
a_{1k}\log|\p^*_n|\, dm$, $A_{2,kn}=\int_\bt a_{2k}\log|\p^*_n|\,
dm$, and recall that $|\psi_n|=1$ a.e.~on $\bt$. Hence, the
function under the exponent sign in \eqref{e7} is purely
imaginary a.e.~on $\bt$. This implies the properties of
$\{A_{0n},A_{1,kn},A_{2,kn}\}$ stated in the lemma, and the proof
of $ii)$ is completed.
\end{proof}

The formulas for coefficients $\{A_{0n},A_{1,kn},A_{2,kn}\}$
in terms of the CMV-represen\-tation can be obtained with the help of the map
$p\mapsto P$ described in \eqref{e101}; the only difference is that $\cc$
should be replaced with its $n\times n$ upper left block $\cc_n$.

\medskip\nt
{\it Proof of Theorem \ref{t2}.} \ We pick a constant $C_1$ in a
way that $0\le C_1p\le 1$ on $\bt$. It is convenient to define
$$
f_n(z)=\exp\lt(\frac12\int_\bt\frac{t+z}{t-z}\, \log\a_n(t)\, dm\rt),\quad
f(z)=\exp\lt(\frac12\int_\bt\frac{t+z}{t-z}\, \log\a(t)\, dm\rt)
$$
with $\a_n(t)=(|\p^*_n(t)|^{-2})^{C_1p(t)}, \a(t)={\s'_{ac}(t)}^{C_1p(t)}$, and
$z\in\bd$.
Obviously, it is enough to show that
$\lim_{n\to\infty}\log f_n(z)=\log f(z)$.
Recalling $\int_\bt |\p^*_n|^{-2}\, dm=1$ \cite{si1}, Theorem 1.7.8, we have
\begin{eqnarray*}
\int_\bt |f_n|^2\, dm&=&\int_\bt \lt(\frac1{|\p^*_n|^2}\rt)^{C_1p}\, dm
=\int_{|\p^*_n|^{-2}\le 1} \lt(\frac1{|\p^*_n|^2}\rt)^{C_1p}\, dm\\
&+&\int_{|\p^*_n|^{-2}\ge 1} \lt(\frac1{|\p^*_n|^2}\rt)^{C_1p}\, dm\le 2
\end{eqnarray*}
It follows similarly that $\int_\bt|f|^2\, dm<\infty$. So, the
functions $f_n,f$ are outer and $\{f_n\}$ is uniformly bounded in
$H^2(\bd)$.
 A ball in $H^2(\bd)$
is weakly compact, and weak convergence  implies the pointwise
convergence on $\bd$. Consequently, there is a subsequence
$\{f_{n_k}\}$ of $\{f_n\}$ that converges to a function $f_0\in
H^2(\bd)$ in $\bd$.

We now prove that $f_0=f$. Indeed, for a $z\in\bd$
$$
\limsup_n \frac12\int_\bt\re\frac{t+z}{t-z}\, p(t)\log\frac1{|\p^*_n(t)|^2}\,
dm(t)\le
\frac12\int_\bt\re\frac{t+z}{t-z}\, p(t)\log\s'_{ac}(t)\, dm(t)
$$
Here, we kept in mind that the measures $|\p^*_n|^{-2}dm$ tend
weakly to $\s$ and the above expressions are semicontinuous with respect to
this
type of convergence \cite[Sect.~5]{ks1}, \cite[Ch.~2]{si1}.  This implies that
$|f_0(z)|\le|f(z)|$ for all $z\in\bd$. We also observe that
$$
\log f_n(0)=\frac12\, \int_\bt(C_1p)\log\frac1{|\p^*_n|^2}\,
dm=\frac12C_1\ti\pps(\cc_n)
$$
where $\ti\pps$ is an expression from \eqref{e4} and $\cc_n$ is the truncated
CMV-matrix. Identity \eqref{e41} from Theorem \ref{t1} reads as
$\log f(0)=(1/2)C_1\ti\pps(\cc)$. In particular, we have
\begin{equation}\label{e8}
\lim_{n\to\infty}\ti\pps(\cc_n)=\ti\pps(\cc)
\end{equation}
which is equivalent to
$$
\log f_0(0)=\lim_{k\to\infty} \log f_{n_k}(0)=\log f(0)
$$
Since the function $f$ is outer, $|f_0|\leq |f|$ and
$|f_0(0)|=|f(0)|$, the usual multiplicative representation of the
functions from $H^2(\bd)$ imply $f=f_0$ on $\bd$. Thus,
the sequence $\{f_n\}$ itself converges to the function $f$, and
the theorem is proved.\epr
\begin{remark}\label{rk1}
Since the function $\eta$ in \eqref{e4} is non-positive,
the convergence in \eqref{e8} is monotone, and $f_{n+1}(0)\le f_n(0)$.
\end{remark}

\section{Asymptotics of orthogonal polynomials in $L^2$-sense}
\label{s3}

For any $\ep>0$, let $B_\ep[\z]=\{z:|z-\z|\le\ep\}$. Furthermore,
let $\oo_\ep=\bd\bsl(\cup_k B_\ep[\z_k])$, $I_{k,\ep}=\bt\cap
B_\ep[\z_k]$, and $A_\ep=\cup_k I_{k,\ep}$. We need several lemmas
to prove the main theorem of this section.
\begin{lemma}\label{l31}
Let $\s\in\ps$. Then, for a finite union of intervals $E\subset \bt$
$$
\limsup_n \int_E p|\log (|\p^*_n|^2\s'_{ac})|\, dm<\infty
$$
\end{lemma}

\begin{proof}
We start by proving that
$$
\limsup_n \int_E p\log^+(|\p^*_n|^2\s'_{ac})\, dm<\infty
$$
Indeed, by $\log^+ x\le x, x>0$ we get
\begin{eqnarray*}
\int_E p\log^+(|\p^*_n|^2\s'_{ac})\, dm&\le&
C\int_E \log^+(|\p^*_n|^2\s'_{ac})\, dm\le
C\int_E |\p_n|^2\s'_{ac}\, dm\\
&\le&C\int_\bt |\p_n|^2\, d\s=C
\end{eqnarray*}
To show  that
$$
\limsup_n \int_E p\log^-(|\p^*_n|^2\s'_{ac})\, dm<\infty
$$
it suffices to know
$$
\liminf_n \int_E p\log(|\p^*_n|^2\s'_{ac})\, dm>-\infty
$$
We have that the measures $\{|\p^*_n|^{-2}\, dm\}$ tend weakly to
$d\s$, and by the semicontinuity of the entropy
\cite[Sect.~5]{ks1}, \cite[Ch.~2]{si1}
$$
\limsup_n\int_E p\log\frac1{|\p^*_n|^2}\, dm\le\int_E p\log\s'_{ac}\, dm
$$
Consequently,
$$
\liminf_n \int_E p\log(|\p^*_n|^2\s'_{ac})\, dm\ge 0
$$
The lemma is proved.
\end{proof}

\begin{lemma}\label{l4}
Let $\s\in\ps$ and
\begin{equation}\label{e9}
\xi_n(z)=\ti D(z)\ti\p^*_n(z)=\exp\lt(\frac12\int_\bt
K(t,z)\log(|\p^*_n|^2\s'_{ac})\, dm\rt)
\end{equation}
Then, for $z\in\oo_{2\ep}$
$$
|\xi_n(z)|\le\frac{C_\ep}{\sqrt{1-|z|}}
$$
where the constant $C_\ep$ does not depend on $n$.
\end{lemma}

\begin{proof} We get $\xi_n=f'_nf''_n$ with
\begin{eqnarray*}
f'_n(z)&=&\exp\lt(\frac12\int_{A_\ep} K(t,z)\log(|\p^*_n|^2\s'_{ac})\,
dm\rt),\\
f''_n(z)&=&\exp\lt(\frac12\int_{\bt\bsl A_\ep}
K(t,z)\log(|\p^*_n|^2\s'_{ac})\, dm\rt)
\end{eqnarray*}
It is plain that, for $t\in A_\ep, z\in\oo_{2\ep}$, the expressions
$|(t+z)/(t-z)|, 1/|q(z)|$ are bounded by constants depending on $\ep$.
Lemma \ref{l31} shows that
$$
\limsup_n \int_{A_\ep} p |\log(|\p^*_n|^2\s'_{ac})|\, dm<\infty
$$
and, therefore, $|f'_n(z)|\le C$ for $z\in\oo_{2\ep}$. Passing to
$f''_n$, we represent it as
\begin{eqnarray}
f''_n(z)&=&\exp\lt(\frac12\int_\bt K(t,z)\log\b_n(t)\rt)\nonumber\\
&=&\exp\lt(\frac12\int_\bt
\frac{t+z}{t-z}\lt(\frac{q(t)}{q(z)}-1\rt)\log\b_n(t)\rt)
\exp\lt(\frac12\int_\bt \frac{t+z}{t-z}\log\b_n(t)\rt)\label{e10}\\
&=&g'_n(z)g''_n(z) \nonumber
\end{eqnarray}
where
$$
\b_n(t)=
\lt\{
\begin{array}{lcl}
1,& \ &t\in A_\ep\\
|\p^*_n|^2\s'_{ac},&& t\in\bt\bsl A_\ep
\end{array}
\right.
$$
Once again, Lemma \ref{l31} implies that
$$
\limsup_n\int_\bt p|\log\b_n|\, dm<\infty
$$
Since $0<c\le p(t)\le C$ for $t\in \bt\bsl A_\ep$, we get
$$
\limsup_n \int_\bt |\log\b_n(t)|\, dm(t)<\infty
$$
Furthermore,
$$
\lt|\frac{t+z}{t-z}\,\frac{q(t)-q(z)}{q(z)}\rt|\le C
$$
for all $z\in\oo_{2\ep}$, and we obtain that $|g'_n(z)|\le C$.

The functions $g''_n$ lie in the Nevanlinna class and are outer. Moreover,
we have
$$
\int_\bt |g''_n|^2\, dm=\int_\bt\b_n\, dm=\int_{\bt\bsl
A_\ep}|\p^*_n|^2\s'_{ac}\, dm\le 1
$$
so $g''_n\in H^2(\bd)$ and $||g''_n||_2\le 1$. To finish the proof
of the lemma, we invoke a standard argument (the integral Cauchy
formula or properties of the reproducing kernel in $H^2(\bd)$)
$$
|g''_n(z)|=|(g''_n, \frac 1{1-\bar zt})|\le||g''_n||_2\,
\lt|\lt|\frac1{1-\bar zt}\rt|\rt|_2 \le \frac{C}{\sqrt{1-|z|^2}}
$$
\end{proof}

The proof of the following lemma is close in spirit to \cite{ki1},
Lemma 3.2.
\begin{lemma}\label{l5} Let $\s\in\ps$. Then
$$
\lim_{n\to\infty} \int_{I'} |\ti D\ti\p^*_n-1|^2\, dm=0
$$
where $I'$ is any closed arc on $\bt$ which does not contain any
point $\{\z_k\}$.
\end{lemma}

\begin{proof}
We fix any closed arc $I$ which does not contain any $\{\z_k\}$
and such that $I'\subset I$. As before, $\xi_n=\ti D\ti\p^*_n$
\begin{figure}[h]
\begin{center}
\includegraphics[height=7cm]{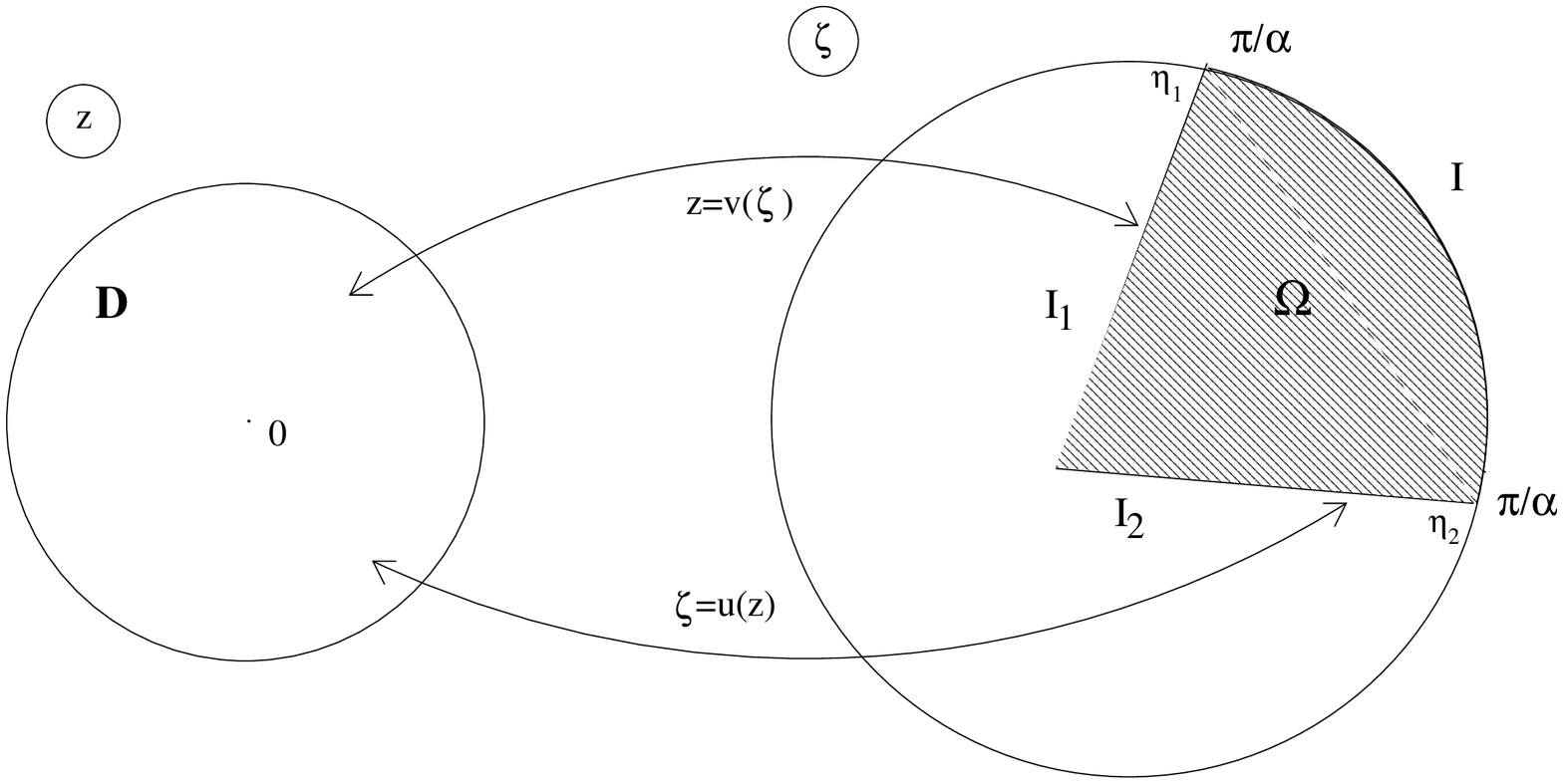}
\end{center}
\caption{}\label{f1}
\end{figure}

Let $\oo$ be the shaded domain on Figure \ref{f1}. Let also
$u:\bd\to\oo$ and $v:\oo\to\bd$ be mutually inverse conformal maps
of the domains, that is, $u(v(\z))=\z$ and $v(u(z))=z$ for
$z\in\bd, \z\in \oo$. We set $\pt\oo$ to be the boundary of $\oo$,
$\pt\oo=I\cup I_1\cup I_2$, where $I$ is the arc on $\bt$ and
$I_1,I_2$ are the straight segments, see the Figure. The angles
between $I, I_1$ and $I_2$ are $\pi/\a, \a>1$. Furthermore, let
$\z_0=u(0)\in\oo$ and $\eta_1,\eta_2$ be the ``corners'' of $\oo$.
It is plain that, for $i=1,2$
\begin{itemize}
\item[\it i)\ \ ] there are constants $c,C>0$ such that
$$
c|\z-\eta_i|^\a\le|v(\z)-v(\eta_i)|\le C|\z-\eta_i|^\a
$$
for $\z\in B_\d(\eta_i)\cap\oo$ and $\d>0$ small enough. \item[\it
ii)\ ] Consequently,
$$
c|\z-\eta_i|^{\a-1}\le|v'(\z)|\le C|\z-\eta_i|^{\a-1}
$$
for these $\z$.
\item[\it iii)]  Obviously,
\begin{itemize}
\item[]\ $c(1-|\z|)^{\a-1}\le|v'(\z)|\le C(1-|\z|)^{\a-1}$ for
$\z\in I_1\cup I_2$, \item[]\ $c|\z-\eta_i|^{\a-1}\le|v'(\z)|\le
C|\z-\eta_i|^{\a-1}$ for $\z\in I$.
\end{itemize}
\end{itemize}
Furthermore, we have
$$
\int_{\pt\oo} |\xi_n(\z)-1|^2|v'(\z)|\, |d\z|=
\int_{\pt\oo}\lt(|\xi_n(\z)|^2-2\re \xi_n(\z)+1\rt)|v'(\z)|\, |d\z|
$$
We start with the second term on the right-hand side
\begin{eqnarray*}
\int_{\pt\oo} 2\re\xi_n(\z)|v'(\z)|\, |d\z|&=&2\re\int_\bt\xi_n(z)\, |dz|\\
&=&4\pi\re\xi_n(u(0))=4\pi\re\xi_n(\z_0)
\end{eqnarray*}
where $\xi_n(z)=\xi_n(u(z))$ and $|dz|=2\pi dm(z)=d\th,\ z=e^{i\th}$. The last
expression in the displayed formula tends to $4\pi$
by Theorem \ref{t2}. Furthermore,
$$
\int_{\pt\oo}|v'|\, |d\z|=\int_\bt |dz|=2\pi
$$
and it remains to show that
\begin{equation}\label{e13}
\lim_{n\to\infty} \int_{\pt\oo}|\xi_n|^2|v'|\, |d\z|\le 2\pi
\end{equation}

We split the last integral into  two integrals over $I$ and
$I_1\cup I_2$, respectively. Then we obtain
\begin{eqnarray*}
\int_I |\xi_n|^2|v'|\, |d\z|&=&\int_I|\p^*_n|^2\s'_{ac}|v'|\,
|d\z|
\le 2\pi\int_I|\p^*_n|^2|v'|\, d\s
\end{eqnarray*}
and the last quantity tends to $\int_I |v'(\z)|\, |d\z|$ by
\eqref{e104}.

We now turn to the integral over $I_1\cup I_2$. Take any $\ep>0$
and freeze it. For any $\d>0$ (its choice will be made precise
later)
$$
\int_{I_1\cup I_2} |\xi_n|^2|v'|\, |d\z|=\int_{I_1\cup I_2,
|\z|\ge1-\d}\ldots+\int_{I_1\cup I_2, |\z|<1-\d}\ldots
$$
and we get for the first integral
\begin{eqnarray*}
\int_{I_1\cup I_2, |\z|\ge1-\d} |\xi_n(\z)|^2|v'(\z)|\,
|d\z|&\le&C\int^\d_0\frac1s\, s^{\a-1}ds=C\int^\d_0 s^{\a-2}\, ds\\
&=&C\d^{\a-1}
\end{eqnarray*}
Above, we used that $\a>1$, a bound from $iii)$ and the inequality
proved in Lemma \ref{l4}. We pick $\d$ small enough to satisfy
$C\d^{\a-1}<\ep$.

Making $\d>0$ smaller, if necessary, we can guarantee that
$$
\lt|\int_{I_1\cup I_2, |\z|\ge 1-\d}|v'|\, |d\z|\rt|<\ep
$$
Then, since $\xi_n$ tends to $1$ uniformly for $|\z|<1-\d$, we take $n$ big
enough to have
$$
\lt|\int_{I_1\cup I_2, |\z|<1-\d}|\xi_n|^2|v'|\, |d\z|-
\int_{I_1\cup I_2, |\z|<1-\d}|v'|\, |d\z|\rt|<\ep
$$
Summing up the inequalities written above, we see that for a large
$n$
$$
\lt|\int_{I_1\cup I_2}|\xi_n|^2|v'|\, |d\z|-
\int_{I_1\cup I_2}|v'|\, |d\z|\rt|<C\ep
$$
which shows
$$
\lim_{n\to\infty} \int_{I_1\cup I_2}|\xi_n|^2|v'|\, |d\z|=\int_{I_1\cup
I_2}|v'|\, |d\z|
$$
So, relation \eqref{e13} is proved. Thus, we obtain
\begin{eqnarray*}
\lim_{n\to\infty}\int_I|\xi_n(\z)-1|^2|v'(\z)|\, |d\z|&\le&\lim_{n\to\infty}
\int_{\pt\oo}|\xi_n(\z)-1|^2|v'(\z)|\, |d\z|\\
&\le&2\pi\lim_{n\to\infty} 2\re(1-\xi_n(\z_0))=0
\end{eqnarray*}
and the lemma is proved for any closed arc $I'\subset I$.
\end{proof}

\begin{remark}\label{rk2}
\hfill
\begin{itemize}
\item[\it i)\ \ ] The lemma also holds for a finite union $A=\cup
I_k$, where $I_k$ are closed arcs that do not contain points from
$\{\z_k\}$. \item[\it ii)\ ] For these arcs $I$, we also have
$$
\lim_{n\to\infty} \int_I |\ti D\ti\p^*_n|^2\, dm=m(I)
$$

\item[\it iii)] For $A\subset\bt$ defined in i),
$$
\limsup_{n\to\infty} \int_{\bt\bsl A} |\xi_n|^2\, dm\le m(\bt\bsl A)
$$
and here  $\{\z_k\}$ necessarily lie in $\bt\bsl A$.
\end{itemize}
\end{remark}

To prove $ii)$ notice that $|||\xi_n||_{L^2(I)}-||1||_{L^2(I)}|\le
||\xi_n-1||_{L^2(I)}$, and the latter quantity tends to $0$ as
$n\to\infty$. As for $iii)$, we have
$$
\lim_{n\to\infty}\int_A|\xi_n|^2\, dm=\lim_{n\to\infty} \int_A
|\p_n|^2\s'_{ac}\, dm=m(A)
$$
so
\begin{eqnarray*}
\limsup_n \int_{\bt\bsl A} |\p_n|^2\s'_{ac}\, dm&\le& 1-\liminf_n
\int_A|\p_n|^2\s'_{ac}\, dm\\
&=&1-\lim_{n\to\infty} \int_A |\p_n|^2\s'_{ac}\, dm=1-m(A)=m(\bt\bsl A)
\end{eqnarray*}

\medskip\nt
{\it Proof of Theorem \ref{t3}.} \ The proof immediately follows
from the Lemma \ref{l5} and the above remarks. Indeed, take an
arbitrary $\ep>0$ and fix it. Then, choose $A=\cup I_k$ (see
$iii)$, Remark \ref{rk2}) in a way that $m(\bt\bsl A)<\ep$.  For
$n$ big enough
$$
\int_{\bt\bsl A} |\xi_n-1|^2\, dm<C\ep
$$
On the other hand, by Lemma \ref{l5}
$$
\lim_{n\to\infty} \int_A |\xi_n-1|^2\, dm=0
$$
and the theorem follows.
\epr

\begin{remark}\label{rk3}
We have
$$
\lim_{n\to\infty} \int_\bt |\p_n|^2\, d\s_s=0
$$
for $\s\in\ps$.
\end{remark}
This is obvious, since $\lim_{n\to\infty} \int_\bt |\p_n|^2\s'_{ac}\, dm=1$ and
$||\p_n||^2_\s=1$.

\section{Modified wave operators and a variational principle}
\label{s4}

\nt {\it Proof of Theorem \ref{t4}.} \ We mainly follow \cite[Ch.~10]{si1}. 
Let us compute $\cf^{-1}\ti\oo_+\cf_0$ on the
vectors of the form $\{z^l\}_{l\in\bz}$; the reasoning for
$\cf^{-1}\ti\oo_-\cf_0$ is similar. Notice that
$A_{0n},A_{2,kn}\in i\br, A_{1,kn}\in\br$ and so the operator
$e^{W(\cc,2n)}$ is unitary. Let $J=\cf^{-1}\cf_0$. Recalling
\eqref{e7}, we get
\begin{eqnarray*}
\cf^{-1}\ti\oo_+\cf_0 z^l&=&\lim_{n\to+\infty}\cf^{-1}e^{W(\cc,2n)}\cc^n\cf\,
(\cf^{-1}\cf_0)\, \cf^{-1}_0\cc^{-n}_0\cf_0 z^l\\
&=&\lim_{n\to+\infty} e^{W(z,2n)}z^nJz^{-n}z^l=
\lim_{n\to+\infty} e^{W(z,2(n+l))}z^{n+l}Jz^{-n}\\
&=&z^l\lim_{n\to+\infty} \psi_{2(n+l)}(z)z^nJz^{-n}
\end{eqnarray*}
and, of course, all limits are to be understood in
$L^2(\s)$-sense. We can assume $n\in\bz_+$ without loss of
generality. Then, $\cf_0 z^{-n}=\cf_0\chi^{(0)}_{2n}=e_{2n}$ and
$\cf^{-1}e_{2n}=z^{-n}\p^*_{2n}(z)$. So
\begin{eqnarray}
\lim_{n\to+\infty} \psi_{2(n+l)}(z)z^nJz^{-n}&=&\lim_{n\to+\infty}
\psi_{2(n+l)}(z)z^{n-n}\p^*_{2n}(z)\nonumber\\
&=&\lim_{n\to+\infty}\psi_{2n}\p^*_{2n}+\lim_{n\to+\infty}
(\psi_{2(n+l)}-\psi_{2n})\p^*_{2n}\label{e131}
\end{eqnarray}
We will prove a little later that
\begin{equation}\label{e14}
\lim_{n\to+\infty} (\psi_{2(n+l)}-\psi_{2n})\p^*_{2n}=0
\end{equation}
in $L^2(\s)$-sense. As for the first term in \eqref{e131}, we have
$$
\lim_{n\to+\infty} \int_\bt |\ti D\psi_n\p^*_n-1|^2\, dm=0
$$
by Theorem \ref{t3}. This is the same as
$$
\lim_{n\to+\infty} \int_\bt |\psi_n\p^*_n-\frac 1{\ti D}|^2\, d\s_{ac}=0
$$
or, together with $\lim_{n\to +\infty}\int_\bt |\p^*_n|^2\, d\s_s=0$ (see Remark
\ref{rk3})$$
\lim_{n\to+\infty} \psi_n\p^*_n=\frac 1{\ti D}\chi_{E_{ac}}
$$
in  $L^2(\s)$-sense, which is exactly the first relation in
\eqref{e141}. Above, $E_{ac}=\bt\bsl\mathrm{supp}\, \s_s$. Let us
prove relation \eqref{e14}. We have
$$
||(\psi_{2(n+l)}-\psi_{2n})\p^*_{2n}||^2_\s=\int_\bt
|\psi_{2n,2(n+l)}-1|^2|\p^*_n|^2\, d\s
$$
where
\begin{eqnarray*}
\psi_{2n,2(n+l)}(z)&=&\exp\Bigg((A_{0\,2(n+l)}-A_{0\,2n})
+\sum_k\Bigg((A_{1,k\,
2(n+l)}-A_{1,k\,2n})\frac{z+\z_k}{z-\z_k}\\
&+&(A_{2,k\,2(n+l)}-A_{2,k\,2n})\lt\{\frac{z+\z_k}{z-\z_k}\rt\}^2\Bigg)\Bigg)
\end{eqnarray*}
Since, by \eqref{e103}, $\lim_{k\to\infty} \a_k=0$ and the coefficients
$A_{0n},A_{1,kn},A_{2,kn}$
depend on a finite number of $\a_k$ only, we have that the expressions in the
small round brackets above tend to zero as $n\to\infty$. Once again, take an
arbitrary $\ep>0$ and fix it. Then, we choose arcs $I'_k, A'=\cup I'_k$, with
the properties $m(A')<\ep$ and $\{\z_k\}\subset A'$. By Remark \ref{rk3} and
$iii)$, Remark \ref{rk2},$$
\int_{A'} |\psi_{2n,2(n+l)}-1|^2|\p^*_n|^2\, d\s\le 4\int_{A'}|\p^*_n|^2\,
d\s<C\ep
$$
for $n$ big enough. On the other hand, $\psi_{2n,2(n+l)}$ uniformly converges
to
$1$ on $\bt\bsl A'$. Hence,
$$
\lim_{n\to+\infty} \int_{\bt\bsl A'} |\psi_{2n,2(n+l)}-1|^2|\p^*_n|^2\, d\s=0
$$
and \eqref{e14} follows. \epr

\medskip
Below, we resort to the notation from the Introduction (see Theorem
\ref{t5}).

\medskip\nt
{\it Proof of Theorem \ref{t5}.} \ We choose a constant $C_0$ from
the condition $C_0\int_\bt p\, dm=1$ and denote the polynomial
$C_0 p$ by $p_0$. Similarly to Lemma \ref{l1}, we have
$$
p_0(t)=a_0+2\re \sum^N_{j=1} a_jt^j
$$
and $a_0=(p_0,1)=\int_\bt p_0\, dm=1$. For any $g\in\cp_0$ and
$j\ge 1$, we have
$$
\int_\bt \log|g|\, dm=\log g(0),\quad \int_\bt\log|g|\, t^j\, dm=
\frac 1{2j!}\, \ovl{(\log g)^{(j)}(0)}
$$
So
$$
\log\l(g)=\log g(0)+\re\sum^N_{j=1} \frac{\bar a_j}{j!}(\log g)^{(j)}(0)
$$
or, what is the same
$$
\l(g)=g(0)\, \exp\lt(\re\sum^N_{j=1} \frac{\bar a_j}{j!}(\log
g)^{(j)}(0)\rt)
$$
Now, if $g\in \cp'_1$, we have $1=\l(g)=||g||_\s\,\l(f)$ \ with
$f=g/||g||_\s\in\cp_0$ and $||f||_\s=1$. Consequently,
$||g||_\s=\l(f)^{-1}$ for these $g,f$ and the infimums in
\eqref{e15} are indeed equal.

For any $g\in \cp_0, ||g||_\s\le 1$, the Jensen inequality implies
$$
\exp\lt(\int_\bt p_0\log\frac{|g|^2\sigma_{ac}'}{p_0}\ dm\rt)\le
\int_\bt |g|^2\, d\s\le 1
$$
This means precisely that
$$
\exp\lt(\int_\bt
p_0\log\frac{\sigma_{ac}'}{p_0}\,dm\rt)\le\exp\lt(-2\int_\bt
p_0\log |g|\, dm\rt)=\frac1{\l(g)^2}
$$
and the first inequality in \eqref{e15} is proved. To deal with
the rest, recall that the measures $(1/|\p^*_n|^2)\, dm$ converge
weakly to $d\s$, and
\begin{eqnarray*}
\liminf_n \int_\bt p_0\log\frac1{|\p^*_n|^2}\, dm &\le& \limsup_n \int_\bt
p_0\log \frac1{|\p^*_n|^2}\, dm\\
&\le& \int_\bt p_0\log \sigma_{ac}'\, dm
\end{eqnarray*}
by the semicontinuity of entropy \cite[Ch.~2]{si1}. The leftmost
expression above is exactly $-2\log\l(\p^*_n)$, and we complete the
proof of the theorem taking exponents in the last inequality.
\epr

\end{document}